\documentclass[11pt]{smfart}

\usepackage[matrix,arrow]{xy}
\usepackage{amscd,amssymb,amsfonts,amsmath}
\usepackage[francais]{babel}
\usepackage{smfthm}
\usepackage{graphics}
\usepackage{epsfig}

\newcommand\mypagesizel{
\textwidth= 6in
\textheight=8.75in
\voffset-.5in
\hoffset -0.35in
\marginparwidth=56pt
}

\newcommand{\Exc}{\textup{Exc}}

\renewcommand{\Big}{\textup{Big}}
\renewcommand{\div}{\textup{div}}
\newcommand{\Div}{\textup{Div}}

\newcommand{\Pef}{\textup{Pef}}

\newcommand{\NS}{\textup{N}^1}
\newcommand{\N}{\textup{N}}
\newcommand{\Z}{\textup{Z}}
\renewcommand{\H}{\textup{H}}

\newcommand{\NE}{\overline{\textup{NE}}}
\newcommand{\Mob}{\overline{\textup{Mob}}}

\newcommand{\epsi}{\varepsilon}
\renewcommand{\phi}{\varphi}

\newcommand{\map}{\dashrightarrow}
\newcommand{\lra}{\to}

\renewcommand{\le}{\leqslant}
\renewcommand{\ge}{\geqslant}

\newcommand{\D}{\Delta}

\newcommand{\mult}{\textup{mult}}

\newcommand{\bP}{\textup{\textbf{P}}}
\newcommand{\bQ}{\textup{\textbf{Q}}}
\newcommand{\bR}{\textup{\textbf{R}}}

\newcommand{\bZ}{\textup{\textbf{Z}}}

\newcommand{\cO}{\mathcal{O}}

\mypagesizel

\begin{document}

\title[Quelques remarques sur la d\'ecomposition de Zariski divisorielle]{Quelques remarques sur la d\'ecomposition de
Zariski divisorielle sur les vari\'et\'es dont
la premi\`ere
classe de Chern est nulle}

\author{St\'ephane DRUEL}

\address{Institut Fourier\\ UMR 5582
du CNRS\\ Universit\'e Joseph Fourier, BP 74\\ 38402 Saint Martin d'H\`eres,
France.}

\email{druel@ujf-grenoble.fr}

\urladdr{http://www-fourier.ujf-grenoble.fr/~druel/}

\maketitle

\section{Introduction}

Toutes nos vari\'et\'es sont
alg\'ebriques et d\'efinies sur le corps des nombres complexes. Soit $X$ une vari\'et\'e projective, lisse et connexe.
On dit qu'un diviseur effectif $D$ sur $X$ a une d\'ecomposition de Zariski s'il
existe des
$\bQ$-diviseurs $P(D)$ et $N(D)$, respectivement num\'eriquement effectif et effectif, tels que 
$$D=P(D)+N(D)$$ 
et tels que l'inclusion
$$\H^0(X,\llcorner mP(D)\lrcorner)\hookrightarrow \H^0(X,\llcorner mD\lrcorner)$$ 
soit 
bijective pour tout entier $m\ge 0$. Zariski \'etablit dans \cite{zariski62} l'existence d'une telle d\'ecomposition
lorsque $X$ est une
surface mais Nakayama montre dans \cite{nakayama04}
qu'en dimension sup\'erieure une
telle d\'ecomposition n'existe
pas toujours m\^eme si on autorise une modification de $X$. 

On a quand m\^eme, \`a la suite des travaux de Nakayama (\cite{nakayama04}) et
Boucksom (\cite{boucksom04}), une d\'ecomposition o\`u la partie positive $P(D)$ est \og num\'eriquement effective en
codimension 1\fg~et
la partie n\'egative est \og rigide\fg~; 
$N(D)$ et $P(D)$ sont \textit{a priori} des $\bR$-diviseurs.

On \'etablit dans cette note quelques propri\'et\'es de cette d\'ecomposition
lorsque la premi\`ere classe de Chern de $X$ est nulle. 

\begin{theo}\label{theo:intro1}
Soit $X$ une vari\'et\'e projective, lisse et connexe avec $c_1(X)=0$.
\begin{enumerate}
\item  Si $D$ un $\bR$-diviseur effectif alors
il existe une vari\'et\'e projective, irr\'eductible et normale $X'$ et une
application birationnelle 
$\phi:X\map X'$ qui contracte les
composantes irr\'eductibles de $N(D)$.
\item Si $D$ est $\bQ$-diviseur grand alors $N(D)$ est un $\bQ$-diviseur.
\end{enumerate}
\end{theo}

On peut bien s\^ur d\'eduire la premi\`ere partie de l'\'enonc\'e du crit\`ere d'Artin (voir \cite{artin62})
si
$\dim(X)=2$.

Boucksom \'etablit dans \cite{boucksom04} la rationalit\'e de la d\'ecomposition de Zariski divisorielle si
$\dim(X)=2$ ou encore si $X$ est une vari\'et\'e hyperk\"ahl\'erienne~: elle se
d\'eduit de l'orthogonalit\'e de cette d\'ecomposition relativement \`a la forme d'intersection dans le premier cas et
\`a la forme de Beauville-Bogomolov dans le second. 

On montre en fait, gr\^ace aux travaux de Birkar, Cascini, Hacon et M$^{\rm c}$Kernan (\cite{BCHM06} et \cite{hm08}),
l'\'enonc\'e plus g\'en\'eral suivant.

\begin{theo}\label{theo:intro2}
Soit $(X,\D)$ une paire klt.
\begin{enumerate}
\item Si $K_X+\D$ est un $\bR$-diviseur pseudo-effectif alors il existe une vari\'et\'e projective, irr\'eductible et
normale $X'$ et une application birationnelle 
$\phi:X\map X'$ qui contracte les
composantes irr\'eductibles de $N(K_X+\D)$.
\item Si $K_X+\D$ est un $\bQ$-diviseur grand alors $N(K_X+\D)$ est un $\bQ$-diviseur.
\end{enumerate}
\end{theo}

On \'etudie dans la derni\`ere partie de ce texte les diviseurs exceptionnels sur les vari\'et\'es
symplectiques holomorphes.
Boucksom a montr\'e que tout diviseur premier exceptionnel sur une 
vari\'et\'e hyperk\"ahl\'erienne est unir\'egl\'e (voir \cite[Proposition 4.7]{boucksom04}). On sait que 
cette propri\'et\'e ne caract\'erise pas
les diviseurs exceptionnels~: on peut penser, par
exemple, au cas d'une
surface $K3$
elliptique g\'en\'erale dont les $24$ fibres singuli\`eres sont des courbes rationnelles (irr\'eductibles).
On obtient le r\'esultat
suivant.

\begin{theo}\label{theo:symplectique}
Soient $X$ une vari\'et\'e projective, lisse et connexe et $E$ un diviseur premier sur $X$. On suppose $X$
symplectique.
Alors $E$ est exceptionnel si et seulement si par un point g\'en\'eral
de $E$ passe une courbe rationnelle $\ell\subset E$ avec $E\cdot \ell<0$.
\end{theo}

On d\'emontre au passage l'\'enonc\'e suivant.

\begin{prop}\label{prop:contraction}
Soient $X$ une vari\'et\'e projective, lisse et connexe et $E$ un diviseur premier sur $X$. On suppose $X$
symplectique et $E$ exceptionnel. Il existe alors un morphisme birationnel $\phi:X\map X'$ avec $X'$ projective,
lisse, connexe et
symplectique et un morphisme birationnel $c':X'\lra Y'$ avec $Y'$ projective et normale dont le lieu exceptionnel est
le support du diviseur premier $E':=\phi(E)$. 
\end{prop}

\section{Multiplicit\'es asymptotiques et d\'ecomposition de Zariski divisorielle}

On commence ce paragraphe par quelques rappels sur les multiplicit\'es asymptotiques et la d\'ecomposition de Zariski
divisorielle. On \'etudie ensuite quelques propri\'et\'es desdites multiplicit\'es.

\subsection{Quelques notations}Soit $X$ une vari\'et\'e alg\'ebrique projective complexe. On suppose $X$ irr\'eductible
et normale.

L'ensemble des diviseurs de Cartier (resp. $\bQ$-diviseurs de Weil $\bQ$-Cartier, $\bR$-diviseurs de
Weil $\bR$-Cartier) est not\'e
$\Div(X)_{\bZ}$ (resp. $\Div(X)_{\bQ}$, $\Div(X)_{\bR}$).
On rappelle que les diviseurs $D_1$ et $D_2$ de $\Div(X)_{\bQ}$ (resp. $\Div(X)_{\bR}$) sont dits
$\bQ$-lin\'eairement \'equivalents (resp. $\bR$-lin\'eairement \'equivalents)
et on note 
$D_1\sim_{\bQ}D_2$ (resp. $D_1\sim_{\bR}D_2$)
s'il existe des fonctions rationnelles $u_j$ non nulles et $r_j\in\bQ$ (resp. $r_j\in\bR$) pour
$j\in J$ fini tels que $D_1-D_2=\sum_{j\in J} r_j\,\div(u_j)$ o\`u $\div(u_j)$ d\'esigne le diviseur des z\'eros
et p\^oles de $u_j$.

On note $\Z_1(X)_{\bZ}$ le groupe ab\'elien libre engendr\'e par les courbes int\`egres et compl\`etes
contenues dans $X$.
On rappelle que les  diviseurs $D_1$ et $D_2$ de $\Div(X)_{\bQ}$ (resp. $\Div(X)_{\bR}$) sont dits
num\'eriquement \'equivalents et on
note 
$D_1\equiv_{\bQ}D_2$ (resp. $D_1\equiv_{\bR}D_2$) si $D_1\cdot C=D_2\cdot C$ 
pour $C\in\Z_1(X)_{\bZ}$.

On note $\N_1(X)$ (resp. $\NS(X)$) l'espace vectoriel r\'eel $\Z_1(X)\otimes_{\bZ}\bR$ 
(resp. $\Div(X)_{\bR}$)
modulo la relation
d'\'equivalence num\'erique d\'efinie ci-dessus. 

Le c\^one convexe ferm\'e de $\N_1(X)$ engendr\'e par les classes des 1-cycles effectifs de $\Z_1(X)_{\bZ}$ est not\'e
$\NE(X)$.

On note $\Pef(X)$ l'adh\'erence dans $\NS(X)$ du c\^one convexe
engendr\'e par les classes des $\bQ$-diviseurs de Weil
effectifs $\bQ$-Cartier.
Un \'el\'ement $D\in\Div(X)_{\bR}$ est dit pseudo-effectif
si sa classe
dans 
$\NS(X)$ est dans $\Pef(X)$.
L'int\'erieur du c\^one $\Pef(X)$ est not\'e $\Big(X)$. Un \'el\'ement $D\in\Div(X)_{\bR}$ est dit grand (\og
big\fg~en anglais)
si sa
classe dans 
$\NS(X)$ est dans $\Big(X)$.

On rappelle enfin qu'un diviseur $D\in\Div(X)_{\bZ}$ est dit mobile si le syst\`eme lin\'eaire correspondant est non
vide et sans composante fixe. On note $\Mob(X)$ le sous-c\^one convexe
ferm\'e de $\NS(X)$ engendr\'e par les classes de diviseurs mobiles. On dit aussi d'un diviseur $D\in\Div(X)_{\bR}$ avec
$[D]\in\Mob(X)$
qu'il est num\'eriquement effectif en codimension 1 ou encore nef en codimension 1

\subsection{La d\'ecomposition de Zariski divisorielle}
On renvoie le lecteur \`a \cite{nakayama04} (voir \'egalement \cite{boucksom04}) pour plus de d\'etails.

Soient $D\in\Div(X)_{\bR}$ et $P$ un diviseur premier sur $X$. 
On suppose $[D]\in\Big(X)$. On d\'efinit la multiplicit\'e asymptotique du syst\`eme lin\'eaire r\'eel
$$|D|:=\{D'\in\Div(X)_{\bR}\text{ effectif tel que }D'\sim_{\bR}D\}$$ en $P$ par
$$\nu_P(D):=\inf_{D'\in |D|_{\bR}}\mult_P(D')\in\bR.$$

On suppose maintenant $[D]\in\Pef(X)$ et on consid\`ere un $\bR$-diviseur $B\in\Div(X)_{\bR}$ avec
$[B]\in\Big(X)$. On a $[D+\epsi B]\in\Big(X)$
pour tout $\epsi>0$ et on
montre que la limite 
$\lim_{\epsi\to 0^+}\nu_P(D+\epsi B)$
existe et ne d\'epend pas de $B$, on la note $\nu_D(P)$.

On montre que $\Mob(X)$ est l'ensemble
des classes des $\bR$-diviseurs pseudo-effectifs $D\in\Div(X)_{\bR}$ tels que $\nu_P(D)=0$ pour tout diviseur premier
$P$ de $X$.

On montre \'egalement que les nombres $\nu_P(D)$ sont nuls sauf pour un nombre fini de
diviseurs premiers $P$ de $X$. On pose alors
$$N(D):=\sum_{P\text{ diviseur premier}} \nu_P(D) P$$
et 
$$P(D):=D-N(D).$$
Le diviseur (effectif) $N(D)$ ne d\'epend que de la classe de $D$ dans $N^1(X)$ et $P(D)\in\Mob(X)$.

L'application obtenue
\begin{eqnarray*}
 \Pef(X) & \longrightarrow & \Mob(X) \\
 {[D]}  & \longmapsto & [P(D)]
\end{eqnarray*}

\noindent est concave, homog\`ene de degr\'e 1 et continue sur $\Big(X)$.

Si $[D]\in\Big(X)$ alors la d\'ecompostion
de Zariski divisorielle est l'unique d\'ecomposition en somme de deux $\bR$-diviseurs $P(D)$ et $N(D)$ respectivement
nef en
codimension un et effectif telle que l'application 
$$\H^0(X,\llcorner mP(D)\lrcorner)\hookrightarrow\H^0(X,\llcorner mD\lrcorner)$$
soit 
bijective pour tout entier $m\ge 0$ o\`u $\llcorner D\lrcorner=\sum_{i\in I}\llcorner d_i\lrcorner D_i$ si
$D=\sum_{i\in I}d_iD_i$ avec $D_i\neq D_j$ pour $i\neq j$ (voir \cite[Theorem 5.5]{boucksom04}).

On rappelle enfin qu'un diviseur effectif $E\in\Div(X)_{\bR}$ est dit exceptionnel si $N(E)=E$ ou, de
fa\c con
\'equivalente, si les classes
des composantes irr\'eductibles de $E$ sont lin\'eairement ind\'ependantes dans $\N^1(X)$ et le c\^one convexe qu'elles
engendrent ne rencontre pas $\Mob(X)$.

\subsection{Quelques propri\'et\'es des multiplicit\'es asymptotiques} Soit $\phi : X\map X'$ une application
birationnelle de vari\'et\'es projectives normales et soit $D\in\Div(X)_{\bR}$. 
On consid\`ere une r\'esolution $\tilde X$ 
des singularit\'es de $X$ et $X'$ avec $q=\phi\circ p$ 
o\`u $p$ et $q$ sont les morphismes de $\tilde X$ sur $X$ et $X'$ respectivement et  
on suppose
$D':=q_*(p^*(D))\in\Div(X')_{\bR}$. On suppose
$|D|_{\bR}$ non vide. On
obtient alors
une application $\bR$-lin\'eaire
\begin{eqnarray*}
|D|_{\bR}  & \longrightarrow & |D'|_{\bR} \\
 G  & \longmapsto &q_*(p^*(G))
\end{eqnarray*}
qui ne d\'epend pas des choix faits. On remarque que si
$\phi^{-1}$ ne contracte pas de diviseur
alors $q_*\circ p^*=\phi_*$.
On souhaite comparer $N(D)$ et $N(D')$.

Soit $P$ un diviseur premier sur $X$.
On suppose que $\phi$ ne contracte pas $P$ et on pose
$P':=\phi_*P$. On suppose d'abord $[D]\in\Big(X)$. On a donc
$[D']\in\Big(X')$. On a aussi
$$\mult_P(G)=\mult_{P'}(q_*(p^*G))\ge\nu_{P'}(D')$$
pour tout $G\in|D|_{\bR}$, et donc
$$\nu_P(D)\ge\nu_{P'}(D').$$
On suppose maintenant $[D]\in\Pef(X)$.
Soient $A$ un diviseur ample sur $X$ et $A':=\phi_*A$. 
On suppose $A'\in\Div(X)_{\bR}$.
On
a 
$$\nu_P(D+\epsi A)\ge\nu_{P'}(D'+\epsi A')$$ 
pour tout $\epsi >0$ puisque $[D+\epsi A]\in\Big(X)$ et,
en passant \`a la limite, on obtient
$$\nu_P(D)\ge\nu_{P'}(D').$$

L'exemple suivant montre qu'on a pas toujours \'egalit\'e dans l'in\'egalit\'e ci-dessus.

\begin{exem}Soit $X'=\bP^1\times C'$ o\`u $C'$ est une courbe compl\`ete (lisse) et soit 
$D'=\bP^1\times\{c'\}$ o\`u $c'$ est un point quelconque de $C'$.
Soit $x'\in D'$ et soit $\phi:X\lra X'$ l'\'eclatement de $X'$ en $x'$. On note
$D:=(\phi^{-1})_*D'$. On a $\nu_P(D)=1$ et $\nu_{P'}(D')=0$ o\`u
$P=D$ et $P'=D'$.
\end{exem}

On \'etudie
maintenant des transformations birationnelles particuli\`eres.

\begin{defi}
Soient $\phi : X\map X'$ une application birationnelle de vari\'et\'es projectives normales et $D\in\Div(X)_{\bR}$. On
dit que $\phi$ est $D$-n\'egative (resp. $D$-strictement
n\'egative) si $\phi^{-1}$ ne contracte pas de diviseur,
$D':=\phi_*D\in\Div(X)_{\bR}$ et
s'il existe une r\'esolution 
$\tilde X$ des singularit\'es de $X$ et $X'$ avec $q=\phi\circ p$ 
o\`u $p$ et $q$ sont les morphismes de $\tilde X$ sur $X$ et $X'$ respectivement pour laquelle
$E:=p^*(D)-q^*(D')$ est effectif (resp. effectif et son support
contient les transform\'es stricts dans $\tilde X$
des diviseurs premiers sur $X$ contract\'es par $\phi$).
\end{defi}

\begin{rema}\label{rema:negativite1}
On suppose que
$\phi : X\map X'$ est $D$-n\'egative (resp. $D$-strictement
n\'egative). On consid\`ere une r\'esolution 
$\tilde X$ des singularit\'es de $X$ et $X'$ avec $q=\phi\circ p$ 
o\`u $p$ et $q$ sont les morphismes de $\tilde X$ sur $X$ et $X'$. On pose
$E:=p^*(D)-q^*(D')$. On montre facilement que $E$ est effectif (resp. effectif et que
son support
contient les transform\'es stricts dans $\tilde X$
des diviseurs premiers sur $X$ contract\'es par $\phi$).
\end{rema}

\begin{rema}\label{rema:negativite2}On suppose \`a nouveau que
$\phi : X\map X'$ est $D$-n\'egative (resp. $D$-strictement
n\'egative).
On consid\`ere $G\sim_{\bR}D$ et $G':=\phi_*G\sim_{\bR}D'$.
On pose $E_G:=p^*(G)-q^*(G')\sim_{\bR} E$. On a en fait
$E_G= E$ par le lemme de n\'egativit\'e \cite[Lemma 3.39]{km98}). On montre de m\^eme que si $X$ et $X'$ sont 
$\bQ$-factorielles alors la propri\'et\'e consid\'er\'ee ne d\'epend que de $[D]\in\NS(X)$ et pas seulement de la
classe de $D$ modulo l'\'equivalence lin\'eaire.
\end{rema}

\begin{lemm}\label{lemm:multiplicite2}
Soit $\phi : X\map X'$ une application birationnelle de vari\'et\'es projectives normales $\bQ$-factorielles. Soient
$D\in\Div(X)_{\bR}$ et $P$ un diviseur premier sur $X$. On suppose $[D]\in\Pef(X)$.
\begin{enumerate}
\item On suppose que $\phi$ est $D$-n\'egative et ne contracte pas $P$. On a alors
$\nu_P(D)=\nu_{P'}(D')$ o\`u $D':=\phi_*D$ et $P':=\phi_*P$.
\item On suppose maintenant que
$\phi$ est $D$-strictement n\'egative et contracte $P$. 
On a alors $\nu_P(D)>0$.
\end{enumerate}
\end{lemm}

\begin{proof}

On consid\`ere une r\'esolution $\tilde X$ 
des singularit\'es de $X$ et $X'$ avec $q=\phi\circ p$ 
o\`u $p$ et $q$ sont les morphismes de $\tilde X$ sur $X$ et $X'$ respectivement.
On pose $E:=p^*(D)-q^*(D')$. On commence par d\'emontrer la premi\`ere assertion de l'\'enonc\'e.
On suppose d'abord $[D]\in\Big(X)$.
Soit $G'\in|D'|_{\bR}$. On a $q^*(G')+E\sim_{\bR} p^*(D)$ et on a donc
$q^*(G')+E=p^*(p_*(q^*(G')+E))$ par le lemme de n\'egativit\'e (voir \cite[Lemma 3.39]{km98}). On en d\'eduit
que l'application lin\'eaire introduite ci-dessus
$$|D|_{\bR}  \longrightarrow |D'|_{\bR} $$
est surjective et
$$\nu_P(D)=\nu_{P'}(D').$$

On suppose maintenant $[D]\in\Pef(X)$. Soit $A'$ un diviseur tr\`es ample sur $X'$
et $A:=(\phi^{-1})_*A'$. Quitte \`a remplacer $A'$ par un diviseur qui lui est lin\'eairement \'equivalent et $A$ par
le diviseur correspondant, on
peut toujours supposer que
$\phi$ est $A$-n\'egative et donc \'egalement 
$(D+A)$-n\'egative.
On
a 
$$\nu_P(D+\epsi A)\ge\nu_{P'}(D'+\epsi A')$$ 
pour tout $\epsi >0$ puisque $[D+\epsi A]\in\Big(X)$ et, en passant \`a la limite, on obtient
$$\nu_P(D)\ge\nu_{P'}(D').$$

On d\'emontre maintenant la seconde assertion de l'\'enonc\'e. On reprend les notations introduites ci-dessus.
On pose $F:=p^*(A)-q^*(A')$ et on note $\tilde P$ le transform\'e strict de $P$ dans
$\tilde X$. 

On fixe $\epsi>0$. On consid\`ere $G\in|D+\epsi A|$ et on pose $G':=\phi_*G$. 
On sait que $\mult_{\tilde P}(E)>0$ par hypoth\`ese et 
on a $p^*(G)=q^*(G')+E+\epsi F$ par la remarque \ref{rema:negativite2}.

On a donc
$$\mult_P(G)=\mult_{\tilde P}(p^*(G)) \ge \mult _{\tilde P} (E+\epsi E_A)\ge \mult_{\tilde P}(E)$$ 
et, en passant \`a la limite, on obtient
$$\nu_P(D)\ge \mult_{\tilde P}(E)>0.$$
\end{proof}

\section{Le programme des mod\`eles minimaux}

On commence ce paragraphe par quelques rappels sur le programme des mod\`eles minimaux ou encore MMP (\og Minimal Model
Program\fg~en anglais). On donne ensuite les d\'emonstrations des th\'eor\`emes \ref{theo:intro1} et \ref{theo:intro2}.

\subsection{Les singularit\'es de paires}

On renvoie pour ce qui suit au tr\`es joli texte \cite{kollar97}.
On d\'esigne par $K_X$ un diviseur canonique sur $X$.
On rappelle
qu'une paire $(X,\Delta)$ est la donn\'ee d'une vari\'et\'e projective $X$ normale et d'un
$\bR$-diviseur de Weil $\D$
sur $X$ tels que $K_X+\D$ soit $\bR$-Cartier.
Soient $(X,\D)$ une paire et $p:\tilde X\lra X$ une r\'esolution des singularit\'es de $(X,\D)$.
On \'ecrit
$$K_{\tilde X}+\tilde{\D}=p^*(K_X+\D)+\sum_{F}a_F(X,\D)F$$
o\`u la somme porte sur l'ensemble des diviseurs premiers $p$-exceptionnels, $\tilde{\D}$ est le transform\'e
strict de $\D$ dans $\tilde X$ et, si $K_X$ est le diviseur d'une forme diff\'erentielle m\'eromorphe
$\omega_X$ sur $X$, $K_{\tilde X}$ est le diviseur de $\omega_X$ sur $\tilde X$. 
Si $F\subset \tilde X$ est un diviseur premier non $p$-exceptionnel, on d\'efinit $a_F(X,\D)$ comme \'etant
l'oppos\'e du coefficient de $F$ dans $\tilde{\D}$.
On appelle le r\'eel $a_F(X,\D)$ la discr\'epance du diviseur premier $F$
relativement \`a la paire $(X,\D)$.

On dit que la paire $(X,\D)$ est \`a singularit\'es terminales (resp. canoniques) si $\D$ est effectif et pour toute
r\'esolution $p:\tilde X\lra X$ des singularit\'es de $(X,\D)$ et tout diviseur premier $p$-exceptionnel $F$,
on a $a_F(X,\D)>0$ (resp. $a_F(X,\D)\ge 0$). On dit que $X$ est \`a
singularit\'es terminales (resp. canoniques) si 
$K_X\in\Div(X)_{\bQ}$ et
$(X,0)$ est \`a
singularit\'es terminales (resp. canoniques).
On rappelle enfin que la paire $(X,\D)$ est dite 
klt (pour Kawamata log-terminale)
si $\D$ est effectif et pour toute
r\'esolution $p:\tilde X\lra X$ des singularit\'es de $(X,\D)$ et tout diviseur premier $F$ de $\tilde X$,
on a $a_F(X,\D)> -1$.

\subsection{Le MMP dirig\'e}\label{mmp}

On renvoie pour ce paragraphe \`a \cite{kmm87}, \cite{km98} et \cite{BCHM06}.

On rappelle qu'un mod\`ele nef (resp. minimal) $(X',\phi)$
de la paire $(X,\D)$ est la donn\'ee d'une paire 
$(X',\D')$ et
d'une application birationnelle $\phi: X\map X'$ 
telles que
\begin{enumerate}
\item $\D'=\phi_*\D$,
\item $\phi$ soit $(K_X+\D)$-n\'egative (resp. $(K_X+\D)$-strictement n\'egative) et
\item $K_{X'}+\D'$ soit nef.
\end{enumerate}

Le MMP dirig\'e est un MMP o\`u les ar\^etes contract\'ees ne sont pas choisies de fa\c con
arbitraire.
Les donn\'ees sont une paire $(X,\D)$ klt o\`u $X$ est $\bQ$-factorielle 
et un $\bR$-diviseur effectif
$H\in\Div(X)_{\bR}$
tel que $K_X+\D+H$ soit nef et $(X,\D+H)$ klt. On pose $t_0=1$.
Le MMP dirig\'e par $H$ produit (conjecturalement) des paires $(X_i,\D_i)$ avec $X_i$ $\bQ$-factorielle, une
suite 
d\'ecroissante de r\'eels $0\le t_i \le 1$ pour $0\le i\le m$, des applications birationnelles 
$\phi_i:X_i \map X_{i+1}$
pour $0\le i\le m-1$ telles que
$K_{X_i}+\D_i+t_iH_i$ soit
nef et $(X_i,\D_i+t_iH_i)$ klt
o\`u, $\D_{i+1}$ (resp. $H_{i+1}$) est le transform\'e strict de $\D_i$ (resp. $H_i$) dans
$X_{i+1}$ et un objet final $(X',\D')=(X_m,\D_m)$ 
tel
que ou bien $(X',\phi)$
soit un mod\`ele minimal de $(X,\D)$, o\`u l'on a pos\'e $\phi:=\phi_{m-1}\circ\cdots\circ\phi_0$, ou
bien $X'$ a une fibration de Mori. Enfin,
$(X_i,\phi_{i-1}\circ\cdots\circ\phi_0)$ est un mod\`ele nef
de la paire $(X,\D+t_iH)$.

On explique maintenant comment faire.
On suppose les $(X_i,\D_i)$ et $0\le t_i\le 1$ d\'ej\`a construits 
avec $K_{X_i}+\D_i+t_iH_i$ nef et on suppose $K_{X_i}+\D_i$ non nef. 
D'apr\`es \cite[Lemma 2.6]{birkar07b}, il
existe
une ar\^ete $R_i\subset\NE(X_i)$ et un r\'eel $0< t_{i+1}\le t_i$ tels que 
$(K_{X_i}+\D_i)\cdot R_i<0$, $(K_{X_i}+\D_i+t_{i+1}H_i)\cdot R_i=0$ et
$K_{X_i}+\D_i+t_{i+1}H_i$ soit nef. Soit
$c_i:X_i\lra Y_i$ la contraction associ\'ee et supposons par exemple la contraction petite.
Soit $c_i^+:X_i^+ \lra Y_i$ le flip de $c_i$ qui existe d'apr\`es \cite[Corollary 1.4.1]{BCHM06}.
On pose $X_{i+1}:=X_i^+$ et
$\D_{i+1}:=\D_i^+$. On a (voir par exemple \cite[Theorem 3.7]{km98}) 
$K_{X_i}+\D_i+t_{i+1}H_i=c_i^*M_i$ o\`u $M_i\in\Div(Y_i)_{\bR}$. Le
diviseur $M_i$ est nef puisque $K_{X_i}+\D_i+t_{i+1}H_i$ l'est, $K_{X_{i+1}}+\D_{i+1}+t_{i+1}H_{i+1}
=((c_i^+)^{-1}\circ c_i)_* (K_{X_i}+\D_i+t_{i+1}H_i)=(c_i^+)^*M_i$ l'est donc \'egalement.
Le cas des contractions divisorielles est analogue.

Si $t\in[0,t_i\lbrack$ ($i\ge 1$),
$(K_{X_{i-1}}+\D_{i-1}+tH_{i-1})\cdot
R_{i-1} <0$ de sorte que $\phi_{i-1}$ est $(K_{X_{i-1}}+\D_{i-1}+tH_{i-1})$-strictement n\'egative par le
lemme
de n\'egativit\'e (\cite[Lemma 3.39]{km98}). On en tire
facilement que l'application
birationnelle $\phi_{i-1}\circ\cdots\circ\phi_0$ est $(K_X+\D+t_iH)$-n\'egative ou encore
que $(X_i,\phi_{i-1}\circ\cdots\circ\phi_0)$ est un mod\`ele nef de $(X,\D+t_iH)$.

\bigskip

Le lemme suivant sera utile au paragraphe \ref{symplectique}.

\begin{lemm}\label{lemm:flip}
On suppose la contraction $c_i$ petite. Alors
$\Exc(c_i^+)$ est couvert par des courbes rationnelles contract\'ees par $c_i^+$.
\end{lemm}
\begin{proof}
On a par construction que $H_i$ est ample$/Y_i$ de sorte que 
$-H_{i+1}$ est ample$/Y_i$. La paire
$(X_{i+1},\D_{i+1}+(t_{i+1}+\epsi)H_{i+1})$ est klt pour $0<\epsi \ll 1$ et 
$-(K_{X_{i+1}}+\D_{i+1}+(t_{i+1}+\epsi)H_{i+1})=-{c_i^+}^*M_i+\epsi(-H_{i+1})$ est ample$/Y_i$.
On d\'eduit alors le r\'esultat cherch\'e de \cite[Theorem 1]{kawamata91}.
\end{proof}

\subsection{Quelques applications}
On ne sait pas d\'emontrer qu'il n'existe pas de suite infinie de flips. On dispose d'un r\'esultat (beaucoup) plus
faible mais suffisant ici.

\begin{prop}
Soient $(X,\D)$ une paire klt avec $X$ $\bQ$-factorielle et $H$ un $\bQ$-diviseur effectif ample sur $X$ tels que 
$K_X+\D+H$ soit nef et la paire $(X,\D+H)$ klt. On consid\`ere un MMP dirig\'e par $H$ pour la paire $(X,\D)$ et on
suppose qu'il n'aboutit pas. On a alors
$$\lim_{i\to+\infty}t_i=0.$$
\end{prop}

\begin{proof}\label{prop:suiteflip}
La suite $(t_i)_{i\ge 0}$ est d\'ecroissante et minor\'ee donc convergente. On suppose que sa limite
$t_{\infty}$ est $>0$. On sait que $(X_i,\phi_{i-1}\circ\cdots\circ\phi_0)$ est un mod\`ele nef de 
$(X,\D+t_iH)=(X,\D+t_{\infty}H+(t_i-t_{\infty})H)$. Or, d'apr\'es \cite[Theorem E]{BCHM06},
l'ensemble des classes d'isomorphie de mod\`eles nef des paires $(X,\D')$ pour 
$\D'\in \D+t_{\infty}H+[0,1-t_{\infty}]H$ est fini.
Il existe
donc deux entiers $i<j$ tels que l'application rationnelle
$\phi_{j-1}\circ\cdots\circ\phi_i$ induise un isomorphisme de $X_i$ sur $X_j$, ce qui
donne la contradiction cherch\'ee, \`a nouveau par le lemme de n\'egativit\'e (voir \cite[Lemma 3.39]{km98}).
\end{proof}

\begin{theo}\label{theo:modelenefcodim1}
Soient $(X,\D)$ une paire klt avec $X$ $\bQ$-factorielle et $H$ un $\bQ$-diviseur effectif ample sur $X$ tels que
$(X,\D+H)$ soit klt. On
suppose $[K_X+\D]\in\Pef(X)$. On consid\`ere un MMP dirig\'e par $H$ pour la paire $(X,\D)$ et on reprend les
notations introduites au paragraphe \ref{mmp}.
On a $K_{X_i}+\D_{i} \in\Mob(X_i)$ pour tout $i \gg 0$, et
les
diviseurs (premiers) contract\'es
sont les composantes irr\'eductibles du support de 
$N(K_X+\D)$.
\end{theo}

\begin{proof}
On a $[K_X+\D]\in\Pef(X)$ et donc,
ou bien le MMP
aboutit a un mod\`ele minimal de la paire $(X,\D)$, ou bien il existe un entier $i_0$ tel que pour tout $i\ge i_0$ les
$\phi_i$ soient des flips de petites contractions.
On commence par le second cas.

Soit $P$ un diviseur premier sur $X$.  On suppose pour commencer que $P$ n'est contract\'e par aucune des applications
rationnelles $\phi_{i-1}\circ\cdots\circ\phi_0$. On note $P_i$ le transform\'e strict de $P$ dans $X_i$.
Par le lemme \ref{lemm:multiplicite2}, on a $\nu_{P}(K_X+\D+t_iH)=\nu_{P_i}(K_{X_i}+\D_i+t_iH_i)$
et $\nu_{P_i}(K_{X_i}+\D_i+t_iH_i)=0$ puisque $K_{X_i}+\D_i+t_iH_i$ est nef par choix de $t_i$.
On d\'eduit finalement de la proposition \ref{prop:suiteflip} que
$\nu_{P}(K_X+\D)=0$. On a donc que $P$ n'est pas une composante irr\'eductible du support de $N(K_X+\D)$.

Inversement, si $P$ n'est pas une composante irr\'eductible du support de $N(K_X+\D)$
alors $\nu_P(K_X+\D)=0$ et $P$
n'est donc
contract\'e par aucune des applications rationnelles $\phi_{i-1}\circ\cdots\circ\phi_0$ par le lemme
\ref{lemm:multiplicite2}.

On d\'eduit \'egalement la premi\`ere assertion de ce que nous venons d'expliquer et on traite le premier cas avec les
m\^emes arguments.
\end{proof}

\begin{coro}Soient $(X,\D)$ une paire klt avec $X$ $\bQ$-factorielle et $H$ un $\bQ$-diviseur effectif ample sur $X$
tels que $(X,\D+H)$ soit klt.
On
suppose $K_X+\D$ de dimension num\'erique $0$, \textit{i.e.} on suppose $[K_X+\D]\in\Pef(X)$ et
$K_X+\D\equiv N(K_X+\D)$. Alors tout MMP dirig\'e par $H$ pour la paire $(X,\D)$ aboutit \`a un mod\`ele minimal
$(X',\D')$ de la paire $(X,\D)$ avec $K_{X'}+\D'\sim_{\bR} 0$.
\end{coro}

On termine ce paragraphe par le r\'esultat suivant.

\begin{theo}\label{theo:rationalite}
Soit $(X,\D)$ une paire klt o\`u $\D$ est un $\bQ$-diviseur
et $[K_X+\D]\in \Big(X)$. Alors $N(K_X+\D)$ est un $\bQ$-diviseur.
\end{theo}

\begin{proof}
Soit $k$ est un entier non nul tel que $k(K_X+\D)$ soit \`a coefficients entiers.
On sait, d'apr\`es
\cite[Corollary 1.1.2]{BCHM06}, que l'alg\`ebre 
$$R(X,k(K_X+\D)):=\bigoplus_{m\ge 0}H^{0}(X,mk(K_X+\D))$$
est de type fini.
Quitte \`a remplacer $k$ par un multiple entier non nul convenable, on peut supposer que l'alg\`ebre
$R(X,k(K_X+\D))$ est engendr\'e par ses \'el\'ements de degr\'e 1, auquel cas, pour tout 
diviseur premier $P$ sur $X$, on a 
$$\nu_P(k(K_X+\D))\in \bZ.$$
On en d\'eduit que les $\nu_P(K_X+\D)$ sont des nombres
rationnels.
\end{proof}

\begin{proof}[D\'emonstration du th\'eor\`eme \ref{theo:intro1}]On applique les th\'eor\`emes
\ref{theo:modelenefcodim1} et \ref{theo:rationalite}
\`a la paire $(X,\epsi D)$ avec $0<\epsi \ll 1$ assez petit.
\end{proof}
\begin{proof}[D\'emonstration du th\'eor\`eme \ref{theo:intro2}]On peut toujours supposer $X$ $\bQ$-factorielle
d'apr\`es \cite[Corollary 1.4.4]{BCHM06}. On d\'eduit les r\'esultats annonc\'es des th\'eor\`emes
\ref{theo:modelenefcodim1} et \ref{theo:rationalite}.
\end{proof}

\section{Cas des vari\'et\'es symplectiques}\label{symplectique}

On montre
facilement qu'un diviseur $E$ sur une vari\'et\'e (lisse) $X$
est exceptionnel si par un point g\'en\'eral de $E$ passe une courbe $C\subset E$ avec $E\cdot C<0$.
On sait d\'emontrer que cette propri\'et\'e caract\'erise les diviseurs exceptionnels si $\dim(X)=2$. On d\'emontre
ici que c'est encore le cas si $X$ est une vari\'et\'e symplectique holomorphe.

\begin{defi}[{\cite[Definition 1.1]{beauville00}}]
Un germe $X$ de vari\'et\'e analytique complexe est dit \`a singularit\'es symplectiques si $X$ est normal et
s'il existe une 2-forme
symplectique
$\omega$ sur le lieu r\'egulier $X_{reg}$ de $X$
telle que, pour toute r\'esolution $p:\tilde X\to X$ des
singularit\'es de $X$, 
$\omega$ s'\'etende en
une $2$-forme r\'eguli\`ere sur $\tilde X$.
\end{defi}

\begin{proof}[D\'emonstration du th\'eor\`eme \ref{theo:symplectique}] 
On fixe un diviseur premier exceptionnel $E$ sur $X$. Soit $0<\epsi \ll 1$ tel que la paire $(X,\epsi E)$ soit klt. On
a $K_X\sim 0$ et $\epsi E=N(K_X+\epsi E)$. Soit $H$ un $\bQ$-diviseur effectif ample sur $X$ tel que la paire $(X,\epsi
E+H)$
soit
encore klt. On sait par le th\'eor\`eme \ref{theo:modelenefcodim1} que tout MMP pour la paire $(X,\epsi E)$
dirig\'e par $H$ contracte $E$. On reprend les notations du paragraphe \ref{mmp}.
Soit $i_0$ tel que le lieu exceptionnel de $c_{i_0}$ soit $E_{i_0}$. 

On sait que $\phi_i$ est un flop pour tout $i\le i_0-1$. On en d\'eduit que $X_i$ est
\`a singularit\'es terminales puis, d'apr\`es \cite[Corollary 1]{namikawa06}
(voir \'egalement \cite{kaledin01}), que $X_i$ est lisse
et symplectique pour tout $i\le i_0$.

On sait par ailleurs que le morphisme $c_{i_0}$
est semi-petit d'apr\`es \cite[Lemma 2.11]{kaledin06} et on a donc $\dim(c_{i_0}(E_{i_0}))=2\dim(X)-2$.
On sait aussi,
d'apr\`es \cite[Theorem 1]{kawamata91}, que $E_{i_0}$ est couvert par des courbes rationnelles
$(\ell_t)_{t\in T}$ telles
que $(K_{X_{i_0}}+\epsi E_{i_0})\cdot \ell_t<0$ ou encore telles que
$E_{i_0}\cdot \ell_t<0$ pour tout $t\in T$. 

Il suffit, pour terminer la d\'emonstration du th\'eor\`eme, de montrer qu'aucune des intersections
$\Exc(c^+_i)\cap E_{i+1}$ pour $i\le i_0-1$ ne domine 
$c_{i_0}(E_{i_0})$ via l'application natuelle 
$E_i\map E_{i_0}\to c_{i_0}(E_{i_0})$.
On suppose que ce n'est pas le cas et on consid\`ere 
un entier $i\le i_0-1$ tel que $\Exc(c^+_i)\cap E_{i+1}$ domine 
$c_{i_0}(E_{i_0})$. On en d\'eduit que l'une des
composantes irr\'eductibles $W_i$ de $\Exc(c^+_i)$ est de 
dimension $2n-2$,
contenue dans $E_{i+1}$ et domine $c_{i_0}(E_{i_0})$~; l'application rationnelle induite 
$W_i\map c_{i_0}(E_{i_0})$
est donc
g\'en\'eriquement finie.
Or, d'apr\`es le lemme \ref{lemm:flip}, $W_i$ est unir\'egl\'ee.
On en d\'eduit que $c_{i_0}(E_{i_0})$
l'est aussi. On sait enfin que la normalis\'ee de
$c_{i_0}(E_{i_0})$ est \`a singularit\'es symplectiques (voir \cite[Theorem 1.4]{wierzba03}). Le lemme \ref{lemm:can}
donne la contradiction cherch\'ee.
\end{proof}

La premi\`ere partie de l'argument ci-dessus d\'emontre la proposition \ref{prop:contraction}.

\begin{lemm}\label{lemm:can}Soit $X$ une vari\'et\'e projective \`a singularit\'es canoniques. 
Si $K_X\sim_{\bQ}0$
alors $X$ n'est pas
unir\'egl\'ee.
\end{lemm}
\begin{proof}Soit $p:\tilde X\lra X$ une r\'esolution des singularit\'es de $X$.
On \'ecrit
$$K_{\tilde X}=p^*(K_X)+\sum_{F}a_F F$$
o\`u la somme porte sur l'ensemble des diviseurs premiers $p$-exceptionnels et $a_F \ge 0$ par hypoth\`ese. On suppose
$X$ unir\'egl\'ee. On note $(\tilde\ell_t)_{t\in T}$ une famille de courbes rationnelles contenues dans $\tilde X$ dont
les d\'eformations dominent $\tilde X$. On a $-K_{\tilde X}\cdot\tilde\ell \ge 2$ d'apr\`es 
\cite[Lemma II.3.13]{kollar96} et donc $-K_{X}\cdot p_*(\tilde\ell) \ge 2$ pour $t\in T$ g\'en\'eral, une contradiction.
\end{proof}

\begin{rema}On reprend les hypoth\`eses de la proposition \ref{prop:contraction}.
On montre facilement (voir \cite[Theorem 1.4]{wierzba03} et \cite[Theorem 4.1]{solwis04}) que
les fibres g\'en\'erales du morphisme $E'\to c'(E')$ sont ou bien des courbes rationnelles lisses ou
bien reunion de deux courbes rationnelles lisses se coupant tranversalement en un point. On peut donc supposer 
$E\cdot \ell=-2$
dans la conclusion du th\'eor\`eme \ref{theo:symplectique}.
\end{rema}

On \'etudie enfin le feuilletage en courbes sur $E$ induit par la forme symplectique ambiante
(voir par exemple \cite{hwangoguiso07}, \cite{hwangviehweg08} et \cite{sawon08}). 

\begin{defi}Soient $X$ une vari\'et\'e (alg\'ebrique) lisse et $E$ un diviseur premier sur $X$. 
On suppose $X$ symplectique et on consid\`ere une forme symplectique $\omega$ sur $X$. Elle induit une $2$-forme
sur l'ouvert dense $E_{reg}$ des points r\'eguliers de $E$ de rang $\dim(E)-1$ dont le noyau d\'efinit un feuilletage en
courbes
appel\'e feuilletage caract\'eristique sur $E_{reg}$ (ou $E$).
\end{defi}

Le r\'esultat suivant g\'en\'eralise \cite[Lemma 10]{sawon08}. On peut d\'eduire 
l'\'enonc\'e de la proposition \ref{prop:feuilletage} du th\'eor\`eme \ref{theo:symplectique} (ou plus exactement de sa
d\'emonstration). On en donne une d\'emonstration \'el\'ementaire.

\begin{prop}\label{prop:feuilletage}
Soient $X$ une vari\'et\'e symplectique, projective et lisse et $E$ un diviseur irr\'eductible sur $X$. 
Si $E$ est unir\'egl\'e
alors les adh\'erences des feuilles g\'en\'erales du feuilletage caract\'eristique sur $E$
sont des courbes rationnelles, et ce sont les seules courbes rationnelles dont les d\'eformations 
dominent $E$.
\end{prop}

\begin{proof}
On consid\`ere une composante irr\'eductible $T$ du sch\'ema des
morphismes $\textup{Hom}(\bP^1,X)$ telle que $ev(T\times\bP^1)$ rencontre $E$ le long d'une partie dense, o\`u 
$ev~:T\times\bP^1\lra X$ est le morphisme universel~; $X$
\'etant symplectique le morphisme $ev$ n'est pas dominant et on a donc $ev(T\times\bP^1)\subset E$.
On en d\'eduit que pour $([f~:\bP^1\to X],b)\in T\times\bP^1$ g\'en\'eral, le rang de la diff\'erentielle
de $ev$ en $([f],b)$ est $\dim(E)=\dim(X)-1$. On a par ailleurs une d\'ecomposition 
$$f^*T_X\simeq\cO_{\bP^1}(a_1)\oplus\cdots\oplus\cO_{\bP^1}(a_{n})
\oplus\cO_{\bP^1}(-a_n)\oplus\cdots\oplus\cO_{\bP^1}(-a_1)$$
avec $a_1\ge \cdots\ge a_{n}\ge 0$ puisque $X$ est suppos\'ee symplectique et enfin
$$rg (d_{([f],b)}ev)= rg(h^{0}(\bP^1,f^*T_X)\otimes\cO_{\bP^1}\lra f^*T_X)$$
d'apr\`es \cite[Proposition II.3.4]{kollar96}. On a donc
$$a_2=\cdots=a_n=0\quad\text{ et }\quad a_1 \ge\deg(T_{\bP^1})=2.$$

On consid\`ere maintenant une r\'esolution $p:X'\lra X$ des singularit\'es de $(X,E)$, on note $E'$ le
transform\'e strict de $E$ dans $X'$ et $f'~:\bP^1\lra X'$ le relev\'e de $f$ \`a $X'$. 

Le point $[f]$ \'etant suppos\'e g\'en\'eral
dans $T$, les d\'eformations de $f'$ passent par un point g\'en\'eral de $E'$ et ${f'}^*T_{E'}$
est donc nef (voir \cite[Proposition II.3.4]{kollar96}). 
On en d\'eduit que l'image de ${f'}^*T_{E'}$ 
dans $f^*T_X$
par l'application
$${f'}^*T_{E'}\subset {f'}^*T_{X'}\lra {f'}^*\circ p^*T_X=f^*T_X,$$
est contenue dans 
$$f^*(T_X)^{\ge 0}:=\cO_{\bP^1}(a_1)\oplus\cO_{\bP^1}^{\oplus(2n-2)}.$$

On a donc, le morphisme $p$ induisant un isomorphisme au-dessus d'un point g\'en\'eral de
$E$, que pour $x\in f(\bP^1)$ g\'en\'eral, l'espace tangent $T_{E,x}$ \`a $E$ en $x$ s'identifie naturellement
au sous-espace 
$(\cO_{\bP^1}(a_1)\oplus\cO_{\bP^1}^{\oplus(2n-2)})\otimes k(x)$ de 
$f^*T_X\otimes k(x)$
puis que $T_{f(\bP^1),x}$ engendre le noyau de la restriction de $\omega$ \`a $E$ en $x$, o\`u 
$\omega$ est une forme symplectique sur $X$,
ou encore
que $f(\bP^1)$ est l'adh\'erence de la feuille du feuilletage caract\'eristique sur $E$.
\end{proof}


\begin{thebibliography}{BCHM06}

\bibitem[Art62]{artin62}
{\scshape M.~Artin} -- {\og Some numerical criteria for contractability of
  curves on algebraic surfaces\fg}, \emph{Amer. J. Math.} \textbf{84} (1962),
  p.~485--496.

\bibitem[BCHM06]{BCHM06}
{\scshape C.~Birkar, P.~Cascini, C.~Hacon {\normalfont \smfandname}
  J.~M\textsuperscript{c}Kernan} -- {\og Existence of minimal models for
  varieties of log general type\fg}, pr\'epublication \'electronique {\tt
  arXiv:math/0610203}, 2006.

\bibitem[Bea00]{beauville00}
{\scshape A.~Beauville} -- {\og Symplectic singularities\fg}, \emph{Invent.
  Math.} \textbf{139} (2000), no.~3, p.~541--549.

\bibitem[Bir07]{birkar07b}
{\scshape C.~Birkar} -- {\og On existence of log minimal models\fg},
  pr\'epublication \'electronique {\tt arXiv:0706.1792}, 2007.

\bibitem[Bou04]{boucksom04}
{\scshape S.~Boucksom} -- {\og Divisorial {Z}ariski decompositions on compact
  complex manifolds\fg}, \emph{Ann. Sci. \'Ecole Norm. Sup. (4)} \textbf{37}
  (2004), no.~1, p.~45--76.

\bibitem[HM08]{hm08}
{\scshape C.~Hacon {\normalfont \smfandname} J.~M\textsuperscript{c}Kernan} --
  {\og Existence of minimal models for varieties of log general type ii\fg},
  pr\'epublication \'electronique {\tt arXiv:0808.1929}, 2008.

\bibitem[HO07]{hwangoguiso07}
{\scshape J.-M. Hwang {\normalfont \smfandname} K.~Oguiso} -- {\og
  Characteristic foliation on the discriminantal hypersurface of a holomorphic
  lagrangian fibration\fg}, pr\'epublication \'electronique {\tt
  arXiv:0710.2376}, 2007.

\bibitem[HV08]{hwangviehweg08}
{\scshape J.-M. Hwang {\normalfont \smfandname} E.~Viehweg} -- {\og
  Characteristic foliation on a hypersurface of general type in a projective
  symplectic manifold\fg}, pr\'epublication \'electronique {\tt
  arXiv:0812.2714}, 2008.

\bibitem[Kal01]{kaledin01}
{\scshape D.~Kaledin} -- {\og Symplectic resolutions: deformations and
  birational maps\fg}, pr\'epublication \'electronique {\tt
  arXiv:math/0012008}, 2001.

\bibitem[Kal06]{kaledin06}
\bysame , {\og Symplectic singularities from the {P}oisson point of view\fg},
  \emph{J. Reine Angew. Math.} \textbf{600} (2006), p.~135--156.

\bibitem[Kaw91]{kawamata91}
{\scshape Y.~Kawamata} -- {\og On the length of an extremal rational curve\fg},
  \emph{Invent. Math.} \textbf{105} (1991), no.~3, p.~609--611.

\bibitem[KM98]{km98}
{\scshape J.~Koll{\'a}r {\normalfont \smfandname} S.~Mori} -- \emph{Birational
  geometry of algebraic varieties}, Cambridge Tracts in Mathematics, vol. 134,
  Cambridge University Press, Cambridge, 1998, With the collaboration of C. H.
  Clemens and A. Corti, Translated from the 1998 Japanese original.

\bibitem[KMM87]{kmm87}
{\scshape Y.~Kawamata, K.~Matsuda {\normalfont \smfandname} K.~Matsuki} -- {\og
  Introduction to the minimal model problem\fg}, Algebraic geometry, Sendai,
  1985, Adv. Stud. Pure Math., vol.~10, North-Holland, Amsterdam, 1987,
  p.~283--360.

\bibitem[Kol96]{kollar96}
{\scshape J.~Koll{\'a}r} -- \emph{Rational curves on algebraic varieties},
  Ergebnisse der Mathematik und ihrer Grenzgebiete. 3. Folge. A Series of
  Modern Surveys in Mathematics [Results in Mathematics and Related Areas. 3rd
  Series. A Series of Modern Surveys in Mathematics], vol.~32, Springer-Verlag,
  Berlin, 1996.

\bibitem[Kol97]{kollar97}
\bysame , {\og Singularities of pairs\fg}, Algebraic geometry---Santa Cruz
  1995, Proc. Sympos. Pure Math., vol.~62, Amer. Math. Soc., Providence, RI,
  1997, p.~221--287.

\bibitem[Nak04]{nakayama04}
{\scshape N.~Nakayama} -- \emph{Zariski-decomposition and abundance}, MSJ
  Memoirs, vol.~14, Mathematical Society of Japan, Tokyo, 2004.

\bibitem[Nam06]{namikawa06}
{\scshape Y.~Namikawa} -- {\og On deformations of {$\Bbb Q$}-factorial
  symplectic varieties\fg}, \emph{J. Reine Angew. Math.} \textbf{599} (2006),
  p.~97--110.

\bibitem[Saw08]{sawon08}
{\scshape J.~Sawon} -- {\og Foliations on hypersurfaces in holomorphic
  symplectic manifolds\fg}, pr\'epublication \'electronique {\tt
  arXiv:0812.3939}, 2008.

\bibitem[SCW04]{solwis04}
{\scshape L.~E. Sol{\'a}~Conde {\normalfont \smfandname} J.~A. Wi{\'s}niewski}
  -- {\og On manifolds whose tangent bundle is big and 1-ample\fg}, \emph{Proc.
  London Math. Soc. (3)} \textbf{89} (2004), no.~2, p.~273--290.

\bibitem[Wie03]{wierzba03}
{\scshape J.~Wierzba} -- {\og Contractions of symplectic varieties\fg},
  \emph{J. Algebraic Geom.} \textbf{12} (2003), no.~3, p.~507--534.

\bibitem[Zar62]{zariski62}
{\scshape O.~Zariski} -- {\og The theorem of {R}iemann-{R}och for high
  multiples of an effective divisor on an algebraic surface\fg}, \emph{Ann. of
  Math. (2)} \textbf{76} (1962), p.~560--615.

\end{thebibliography}

\providecommand{\bysame}{\leavevmode ---\ }
\providecommand{\og}{``}
\providecommand{\fg}{''}
\providecommand{\smfandname}{et}
\providecommand{\smfedsname}{\'eds.}
\providecommand{\smfedname}{\'ed.}
\providecommand{\smfmastersthesisname}{M\'emoire}
\providecommand{\smfphdthesisname}{Th\`ese}

\end{document}